\documentclass{article}

\usepackage{amssymb, amsmath, amsthm, amsfonts, amscd, epsf, graphicx}
\usepackage[dvips,usenames]{color}
\usepackage[american]{babel}
\usepackage[latin1]{inputenc}
\usepackage[T1]{fontenc}
\usepackage{ae, aecompl}
\usepackage[mathscr]{eucal}

\newtheorem{theorem}{Theorem}

\newtheorem{corollary}{Corollary}
\newtheorem{proposition}{Proposition}

\newtheorem{ex}{Example}
\newtheorem{remark}{Remark}

\def\Z{\mathbb Z}

\def\Q{\mathbb Q}

\def\a{\alpha}

\def\<{\langle}
\def\>{\rangle}

\def\Gal{\mbox {Gal }}

\def\D{\Delta}

\def\Gal{\mbox{Gal }}

\def\deg{\mbox{ deg }}

\def\a{\alpha}

\def\<{\langle}
\def\>{\rangle}

\def\a{\alpha}

\title{Galois groups of prime degree polynomials with nonreal roots}

\author{A. Bialostocki\\
Department of Mathematics\\
University of Idaho\\
Moscow, ID, 83843. \\ \\
and\\ \\
T. Shaska\footnote{Supported by the NSA grant R1-05-0129}\\
Department of Mathematics\\
Oakland University\\
Rochester, MI, 48309-4485.}

\begin{document}

\maketitle

\begin{abstract}
In the process of computing the Galois group of a prime degree polynomial $f(x)$ over $\Q$ we suggest a
preliminary checking for the existence of non-real roots. If $f(x)$ has non-real roots, then combining a 1871
result of Jordan and the classification of transitive groups of prime degree which follows from CFSG we get that
the Galois group of $f(x)$ contains $A_p$ or is one of a short list. Let $f(x)\in \Q[x] $ be an irreducible
polynomial of prime degree $p \geq 5$ and $r=2s$ be the number of non-real roots of $f(x)$. We show that if $s$
satisfies $s \, ( s \log s + 2 \log s + 3) \leq p $ then $Gal (f)= A_p, S_p$.
\end{abstract}

\section{Introduction}
Solving algebraic equations is one of the oldest and most fundamental problems in mathematics. The problem was put
on a firm basis with the contribution of E. Galois and was the main motivation for the development of modern
algebra. It is now considered basic knowledge that a polynomial equation with rational coefficients can be solved
by radicals if and only if its Galois groups is solvable. Hence, the problem of determining Galois groups is
important and has been considered in many areas of mathematics as number theory, group theory, algebraic geometry,
and differential equations.

However, computing Galois groups is still a difficult task. Even with the development of new computer algebra
systems this remains a challenge and can be accomplished only for small degree polynomials. For example, Maple 9
can only handle polynomials of degree $\leq 9$ and Kant up to degree 15. Other computer algebra packages can
handle polynomials whose degree is in the same range.

Let $f(x)=a_0 + a_1 x + \cdots + a_n x^n=0$ be an algebraic equation. Suppose that $a_i= \pm 1$; then Littlewood
and Offord proved that almost all these equations have less than $25(\log n^2)$ real roots; see \cite{LO}. If the
$a_i$ have normal distribution with density $e^{-u^2}/\pi^{{\textstyle\frac 1{2}}}$; then Kac (1944) showed that
the mean value of the number of real roots of $f(x)$ equals
$$
(4/\pi)\int_0^1\frac{[1-n^2[x^2(1-x^2)/(1-x^{2n})]^2]^{{\textstyle\frac 1{2}}}}{1-x^2}\,dx\sim(2/\pi)\log n ,$$
see \cite{Kac1}. If the coefficients are random variables with distribution function $\sigma(u)$ then in a second
paper Kac (1949) gives an expansion for the average number of roots of the above equation in the interval $(a,b)$.
He then considers the special case when $\sigma(u)={\textstyle\frac 1{2}}u$, and deduces from the general formula
that in this case the average number of real roots is asymptotic to $(2\pi)^{-1}\log n$. Consequently, polynomials
in general have plenty of non-real roots.

The existence of non-real roots of a polynomial makes the computation of its Galois group much easier. Computing
the Galois group in this case for polynomials of prime degree $p$ will be the focus of this short note. Checking
whether a polynomial has non-real roots is very efficient since numerical methods can be used. Once the existence
of non-real roots is established then from a theorem of Jordan (1871) it follows that if their number is "small"
enough with respect to the degree $p$ of the polynomial, then the Galois group is $A_p$ or $S_p$. Furthermore,
using the classification of finite simple groups we know nowadays the complete classification of transitive groups
of prime degree. This enables us to provide a complete list of possible Galois groups for every polynomial of
prime degree $p$ which has non-real roots.

In section 2 we briefly describe the existing techniques used to compute the Galois group of a polynomial $f(x)
\in \Q[x]$ of degree $n$. These techniques are mainly based on the Dedekind's theorem and knowledge of the list of
all transitive subgroups of $S_n$ and the cycle structure of their elements. First, a factorization $\mod p$
(i.e., good prime $p$) of $f(x)$ is required to obtain information on the cycle structure of the group elements.
Such information enables us to eliminate groups from the list of transitive subgroups of $S_n$. For as long as the
group is not uniquely determined we repeat the procedure with a different prime $p$. This is a rather expensive
technique since algorithms of factorizing polynomials are not very efficient. Furthermore, many primes $p$ might
be needed in the process.

In section 3 we study polynomials of prime degree $p$ with non-real roots. We describe the transitive groups of
$S_p$ and compute all such groups for $n \leq 30$. Using a theorem of Jordan on permutation groups we show that
for a fixed number of non-real roots $r=2s$ and
$$p \geq  N(r):= \left[ s \, ( s \log s + 2 \log s + 3) \right] $$
the Galois group is $A_p$ or $S_p$. Furthermore, we classify all groups that occur for $p < N(r)$ and provide an
algorithm that computes the Galois group of a polynomial with non-real roots of degree prime $p$. We conclude with
some final remarks in section 4.

\medskip

\noindent \textbf{Notation: } Throughout this paper the ground field $k=\Q$. All polynomials are assumed to be
irreducible over $\Q$. The group $D_n$ denotes the dihedral group of $n$ elements, $M_{11}$ and $M_{23}$ are
respectively Mathew groups of degree 11 and 23. For all other groups we use the GAP notation.

\section{The Galois group of a polynomial}
Let $f(x)\in k[x]$ be a degree $n$ polynomial. A \textit{splitting field} of $f$ is a field extension of $k$ of
the form $k(\a_1, \dots ,\a_n)$ where $f(x)=(x-\a_1)\cdots (x-\a_n)$. Any two splitting fields of $f$ are
isomorphic under an isomorphism trivial on $k$. Thus we normally speak of \textit{the splitting field} of $f$.

If the discriminant $\D_f $ of $f(x)$ is nonzero then $f(x)$ has $n$-distinct roots $\a_1, \dots ,\a_n$ in the
splitting field $E_f$ of $f$. $E_f/k$ is a Galois extension. The \textit{Galois group} of $f$ over $k$, denoted by
$G_k(f)$, is the group $G(E_f/k)$, viewed as a permutation group of the roots $\a_1, \dots ,\a_n$. Thus $G_k (f)$
is a subgroup of $S_n$, determined up to conjugacy by $f$. The following is elementary and we avoid the proof.

\begin{proposition} Let $f(x)\in k[x]$ be a degree $n$ polynomial and $G=G_k(f)$.\\
 (i) Let $H=G\cap A_n$. Then $H=\ G(E_f / k (\sqrt{\D_f}))$. In particular, $G$ is contained in the
alternating group $A_n$ if and only if the discriminant $\D_f$ is a square in $k$.\\
(ii) The irreducible factors of $f$ in $k[x]$ correspond to the orbits of $G$. In particular, $G$ is a transitive
subgroup of $S_n$ if and only if $f$ is irreducible.
\end{proposition}

\begin{remark}
Recall that for a degree $n$ irreducible polynomial $f(x)$ over $k$ with splitting field $E_f$ we have $n \mid
[E_f : k]$.
\end{remark}

\noindent Hence, for a given polynomial $f(x)$ of degree $n$ its Galois group must satisfy:

i) $G$ is isomorphic to a transitive subgroup of $S_n$.

ii) $n$ divides $ |G|$,

iii) $G $ is a subgroup of $A_n$ if and only if $\D_f$ is a square in $\Q$.

\medskip

\noindent These conditions narrow down the possible choices of groups that can be Galois groups of $f(x)$. In
order to determine precisely the group $G$ we need to determine the type of cycles in $G$.

\subsection{Reduction mod p}

The reduction method uses the fact that every polynomial with rational coefficients can be transformed into a
monic polynomial with integer coefficients without changing the splitting field.

Let $f(x) \in \Q[x]$ be given by
$$f(x)= x^n + a_{n-1} x^{n-1} + \cdots + a_1 x + a_0$$
Let $d$ be the common denominator of all coefficients $a_0, \cdots , a_{n-1}$. Then $g(x):= d\cdot f(\frac x d)$
is a monic polynomial with integer coefficients. Clearly the splitting field of $f(x) $ is the same as the
splitting field of $g(x)$. Thus, without loss of generality we can assume that $f(x)$ is a monic polynomial with
integer coefficients. The following theorem gives information on the cycle shape of permutations of $\Gal (f)$.

\begin{theorem}\textbf{(Dedekind)} Let $f(x) \in \Z[x]$ be a monic polynomial such that  $\deg f = n$, $\Gal_\Q (f) = G$,
and $p$ a prime such that $p \nmid \D_f$. If $f_p:=f(x) \mod p \, \, $ factors in $\Z_p [x]$ as a product of
irreducible factors of degree
$$n_1, n_2, n_3, \cdots , n_k,$$
then $G$ contains a permutation of type
$$(n_1)\, (n_2) \, \cdots \, (n_k)$$
\end{theorem}

From the Chebotarev density theorem we know that if $p \to \infty$ then the distribution of factor degrees
approaches the distribution of the cycle shapes in the group. Assuming that we can compute all transitive
subgroups of $S_n$ and their cycle shapes then the above two theorems give a basis of an algorithm to determine
$Gal (f)$. The transitive subgroups of $S_n$ can be computed in GAP for all $n\leq 30$.

\begin{remark}
The above technique doesn't determine the Galois group uniquely in all cases. For example there are two
non-isomorphic degree 8 groups with the same cycle structure. In such cases other methods such as invariants of
groups are used to determine the group uniquely.
\end{remark}

\subsection{Transitive subgroups of $S_n$}
By \emph{degree of a permutation group} $G$ we mean the number of points moved by $G$. The \emph{degree of a
permutation} $\a \in S_n$ is the number of points moved by $\< \a \>$. The \emph{minimal degree} of $G$, denoted
by $m (G)$, is the smallest of degrees of elements $\a \neq 1$ in $G$.

In the case $n$ is a prime $n=p$ then a non-solvable transitive subgroup $G$ of $S_p$ is doubly transitive
(Burnside theorem); see \cite[pg. 431]{Bu} or \cite[Theorem 11.7]{Wi}. Using the classification of simple groups
one gets a complete list of $G$:

\begin{theorem} Let $G$ be a doubly transitive subgroup of $S_p$, for a prime $p$. Then $G$ isomorphic to one of the
following:\\
(i) $A_p, S_p,$\\
(ii) $p = 11, G = L_2(11)$ or $M_{11}$,\\
(iii) $p = 23,\, \, G = M_{23}$,\\
(iv) $p = \frac {(q^k-1)} {(q-1)}$ and $L_k (q) \leq G \leq Aut(L_k(q))$.
\end{theorem}

\begin{proof} One can check Prop.~4.4.1 in \cite{KT} for all these statements.
\end{proof}

Thus, the above theorem classifies all transitive non-solvable subgroups of $S_p$. For the solvable cases we use
GAP to compute them in each case for $\ \leq 29$. In the Table 1, we display the number of transitive subgroups
and the number of non-solvable transitive subgroups of $S_n$ for all $n\leq 30$.
\begin{table}[ht!]     \label{table}
\vspace{2mm}
\begin{center}
\begin{tabular}{||c|c|c||}
\hline  \hline && \\
Deg. & Nr. of trans. groups & Nr. of unsolvable groups \\
 & & \\
\hline & & \\
5 & 5  &  2   \\
6 & 16  &  4   \\
7 & 7  &  3   \\
8 & 50  &  5   \\
9 & 34  &  4   \\
10 & 45  & 21 \\
11 & 8  &  4   \\
12 & 301  &  36   \\
13 & 9  &  3   \\
14 & 63  &  27   \\
15 & 104  &  40 \\
16 & 1954  &  49   \\
17 & 10  &  5   \\
18 & 983  &  91   \\
19 & 8  &  2   \\
20 & 1117  &  358   \\
21 & 164  &  56   \\
22 & 59  &  27   \\
23 & 7  &  3   \\
24 & 25000  &  807   \\
25 & 211  &  79 \\
26 & 96  &  26   \\
27 & 2392  &  64   \\
28 & 1854  &  617   \\
29 & 8  &  2 \\
30 & 5712  &  1896 \\
\hline \hline
\end{tabular}
\vspace{3mm} \caption{The number or transitive and unsolvable groups (up to conjugacy) for degree $\leq 30$}
\end{center}
\end{table}

\section{Polynomials with non-real roots}

Let $f(x)\in \Q[x] $ be an irreducible polynomial of degree $n > 5$. Denote by $r$ the number of non-real roots of
$f(x)$. Since the complex conjugation permutes the roots then $r$ is even, say $r=2s$. By a reordering of the
roots we may assume that if $f(x)$ has $r$ non-real roots then
$$\a:=(1, 2) (3, 4)\cdots (r-1, r) \in Gal (f).$$
Since determining the number of non-real roots can be very fast, we would like to know to what extent the number
of non-real roots of $f(x)$ determines $Gal (f)$. The complex conjugation assures that $m(G) \leq r$. The
existence of $\a$ can narrow down the list of candidates for $Gal (f)$. However, it is unlikely that the group can
be determined only on this information unless $p$ is "large" enough. In this case the number of non-real roots of
$f(x)$ can almost determine the Galois group of $f(x)$, as we will see in the next section. Nevertheless, the test
is worth running for all $p$ since it is very fast and improves the algorithm overall.

\subsection{Polynomials of prime degree}
An approach of computing Galois groups of polynomials or solving for roots is to check whether or not the
polynomial is decomposable. The polynomial decomposition problem can be stated as follows: given a degree $n$
polynomial $f \in k[x]$, determine whether there exist polynomials $f_1, f_2 $ of degree greater than one such
that $f=f_1 \circ f_2=f_1 ( f_2(x))$, and in the affirmative case to compute them. From the classical L\"uroth's
theorem this problem is equivalent to deciding if there exists a proper intermediate field in the finite algebraic
extension $k(f) \subset k(x)$. From the computational point of view, there are several polynomial time algorithms
for decomposing polynomials. The computation of $f_1(x)$ and $f_2 (x)$ only requires $O(n^2)$ arithmetic
operations in the ground field $k$; see for instance \cite{Gu}. A polynomial $f(x)\in F[x] $ is indecomposable
over the subfield $F \subset k$ if and only if $f(x)$ is indecomposable over $k$. There are fast algorithms to
compute the decomposition of polynomials; see \cite{Gu}. However, for indecomposable polynomials we would like to
have better methods of determining the Galois group and possible roots of the polynomial. Hence, polynomials of
prime degree are of interest since they are, of course, indecomposable.

The next theorem determines the Galois group of a prime degree polynomial $f(x)$ with $r$ non-real roots when the
degree of $f(x)$ is large enough with respect to $r$.

\begin{theorem}
Let $f(x)\in \Q[x] $ be an irreducible polynomial of prime degree $p \geq 5$ and $r=2s$ be the number of non-real
roots of $f(x)$. If $s$ satisfies
$$s \, ( s \log s + 2 \log s + 3) \leq p $$
then $Gal (f)= A_p, S_p$.
\end{theorem}

\begin{proof} Since $p$ is prime then every transitive subgroup of $S_p$ is primitive. Let $G$ denote the
Galois group of $f(x)$ and $m (G)$ its minimal degree. By reordering the roots we can assume that
$$(1, 2) (3, 4)\cdots (r-1, r) \in Gal (f).$$
Hence, $m:=m (G) \leq r$. From a theorem of Jordan \cite{Jo} we have that if
$$ \frac {m^2} 4 \log \frac m 2 + m \left( \log \frac m 2 + \frac 3 2 \right) \leq p$$
then $G= A_p$ or $S_p$. Hence, if $$s \, ( s \log s + 2 \log s + 3) \leq p $$ then $Gal (f)= A_p$ or $S_p$.
\end{proof}

\begin{remark} For a modern view of Jordan's theorem and its implications to number theory, Galois
representations, and topology see \cite{Se1}.
\end{remark}

For a fixed $p$ the above bound is not sharp as we will see below. However, the above theorem can be used
successfully if $s$ is fixed. We denote the above bound on $p$ by
$$N(r) := \left[ s \, ( s \log s + 2 \log s + 3) \right] $$
for $r=2s$. Hence, for a fixed number of non-real roots, for $p \geq N(r)$ the Galois group is always $A_p$ or
$S_p$.

\begin{corollary} Let a polynomial of prime degree $p$ have $r$ non-real roots. If one of the following
holds:

(i) $r=4$ and $p > 7$,

(ii) $r=6$ and $p > 13$,

(iii) $r=8$ and $p > 23$,

(iv) $r=10$ and $p > 37$,\\

\noindent then $Gal (f) = A_p$ or $S_p$.
\end{corollary}

\begin{remark}
The above results gives a very quick way of determining the Galois group for polynomials with non-real roots.
Whether or not the discriminant is a complete square can be used to distinguish between $A_p$ and $S_p$.
\end{remark}

If $p < \ N(r) $ then some exceptional cases occur. Next theorem determines these exceptional cases for
polynomials of degree up to 29. The computations were made using GAP.

\begin{theorem}
Let $f(x)\in \Q[x] $ be an irreducible polynomial of prime degree $p > 5$. Let $r$ be the number of complex roots
of $f(x)$. If $r > 0$ then $Gal (f) $ is $A_p$, $S_p$ or one of the groups as in the following Table 2.
\end{theorem}

\begin{tiny}
\begin{table}[ht!]
\caption{Galois groups (other then $A_p, S_p$) of polynomials with non-real roots.}
{ \begin{tabular}{||c|c|c||c|c||}
\hline \hline
 $p$ &  Solv. & Sign. & Nonsol.  & Sign. \\
\hline \hline
7  & $ D_7  $      &  $ (2)^3, (7)  $   & $L(7)$   &     $ (2)^2, (4)(2), (3)^2, (7) $ \\
   & $(7, 4)$   &  $ (2)^3, (3)^2,  (7)  $    &         &   \\
\hline \hline
11 & $ D_{11}$     &  $ (2)^5,  (11)  $    & $L(11) $        & $(2)^4, (3)^3, (5)^2, (2)(6)(3), (11) $\\
   & & & & \\
   & $ (11, 4)$&  $  (2)^5, (5)^2, (10), (11) $    & $ M_{11}$       & $ (2)^4, (2) (6)(3), (2) (8), (3)^3$,  \\
      & & & &  $(4)^2, (5)^2, (11) $\\
\hline \hline
13 & $D_{13}$      &  $ (2)^6, (13)  $    & $L(13)$        & $ (2)^4, (3)^3, (3)^4, (4)^2 (2)^2,$  \\
   & & & & $ (6)(3)(2), (8)(4), (13)  $\\
   & $ (13, 4)$ &  $ (2)^6, (4)^3, (13)  $    &                 & \\
   & $ (13, 5)$ &  $ (2)^6, (3)^4, (6)^2, (13)$    &                & \\
   & $ (13, 6)$&  $ (2)^6, (3)^4, (4)^3$,     &                & \\
   &               &  $(6)^2, (12), (13) $         & & \\
\hline \hline
17 & $D_{17}$      &  $ (2)^8, (17)  $    & $PSL_2 (16)$   & $ (2)^8, (3)^5, (5)^3, (15), (17)  $ \\
   & & & & \\
   & $ (17, 3)$ &  $  (2)^8, (4)^4, (17) $    & $(17, 7)$     & $ (2)^6, (2)^8, (3)^5, (4)^4,  $ \\
   & & & & $ (5)^3, (6)^2 (3), (5)(10)(2), $\\
   & & & & (15), (17) \\
      & & & & \\
   & $ (17, 4)$&  $ (2)^8, (4)^4, (8)^2, (17)  $    &$ (17, 8)$   & $ (2)(5)(10),  (2)(4), (2)(4)^3, $ \\
   & & & & $ (2)^6, (2)^8, (3)(6)^2, (3)^5,  $\\
   & & & & $ (3)^2 (12), (4)^3, (5)^3, (8)^2, $\\
     & & & & (15), (17)  \\
   & $ (17, 5)$&  $  (2)^8, (4)^4, (8)^2, (16), (17)  $    &                & \\
\hline \hline
19 & $D_{19}$      & $ (2)^9, (19)  $     &                & \\
   & $ (19, 4)$& $ (2)^9, (3)^6, (6)^3, (19) $     &  & \\
   & $ (19, 6)$& $ (2)^9, (3)^6, (6)^3, (9)^2,   $     &  & \\
   &               & $ (18), (19)  $ & & \\
\hline \hline
23 & $D_{23}$      & $ (2)^{11}, (23)   $     & $M_{23}$       & $ (2)^8, (2)^2 (4)^4, (2) (7) (14),    $ \\
 & && & $ (2) (4) (8)^2, (2)^2 (3)^2 (6)^2, $ \\
&$(23, 4)$ &$ (2)^{11}, (11)^2, (22), (23)    $& & $(3) (5) (15), (5)^3, (5)^4, $ \\
& && &$(7)^3, (11)^2, (23)   $ \\
\hline \hline
29 & $D_{29}$      &  $(2)^{14}, (29)   $    & & \\
   & $(29, 3) $       &  $ (2)^{14},  (4)^7, (29)  $    & & \\
   & $ (29, 5)$      &  $  (2)^{14},  (7)^4, (14)^2,  (29)    $    &  & \\
   & $ (29. 6)$&  $ (2)^{14}, (4)^7, (7)^4, (14)^2,   $    &  & \\
   & & $ (28), (29)$ & & \\
\hline \hline
\end{tabular}
}
\end{table}

\end{tiny}

\proof The proof is computational and follows from the tables of transitive subgroups of $S_p$. It is easy to
decide which ones of these groups are nonsolvable and compute their cycle types.
\endproof

\begin{remark} We used in  Table 2 notations which we considered standard as $D_n$, $M_{11}$,
and $L(p)$, otherwise we used the GAP notation $(p, i)$ which is the i-th group in the list of transitive groups
of degree $p$. These groups can be generated in GAP simple by typing \verb"TransitiveGroup(n,i);". The group
$M_{23}$ is not realized as a Galois groups over $\Q$.
\end{remark}

Notice that no two groups of Table 2 have the same cycle structure. Hence the Galois group can be determined
uniquely by reduction mod $p$ for all polynomials of prime degree $\leq 29$ with non-real roots.

\begin{ex}
Let $f(x)$ be the polynomial
$$f(x)=x^{11}+5 x^7-4 x^6-20 x^5+4 x^4+20 x^3+1$$
This polynomial is irreducible over $\Q$ and has exactly 10 non-real roots. The reader can easily
check these facts in Maple using the commands:\\
\verb" factor(f(x));" \verb" realroot(f(x));"

\noindent From the above theorem, its Galois group is $A_{11}$ or $S_{11}$. Its discriminant is
$$\Delta_f =  - 59 \cdot  1391212936091429123033$$
which is not a square in $\Q$. Hence the Galois group of $f(x)$ is $S_{11}$. Maple can not compute the Galois
group of this polynomial since its degree is $ > 8$.
\end{ex}

Combining the above results we have the following algorithm for computing the Galois group of prime degree
polynomials with non-real roots. \verb"D(f)" denotes the discriminant of $f(x)$ and \verb"A_p", \verb"S_p" the
alternating and symmetric group of $p$ letters. Note that in the case $p< N(r)$ the we know that a permutation of
the type $(2)^{\frac r 2}$ is in the group. Hence, the list of transitive subgroups is much shorter than in
general. This information was obtained by computing the number of real roots other then by some factorization
modulo $p$. Thus, even in this case the algorithm is improved.

\bigskip

\vskip3truept\hrule \vspace*{1.8 mm}

\noindent {\sc Algorithm:} Computing the Galois group of prime degree polynomials with non-real roots.

\vskip3truept\hrule\vskip4truept

\medskip

\noindent {\bf Input:} {\rm An irreducible monic polynomial $f(x) \in \Q [x]$ of prime degree $p$.}

\noindent {\bf Output:} {\rm Galois group $Gal (f)$ of $f(x)$ over $\Q$ }

\medskip

\begin{verbatim}
begin

r :=NumberOfRealRoots(f(x));

If p > N(r) then
   if  D(f) is a complete square then
       Gal (f)=A_p; else Gal (f) = S_p;
   endif;
else ReductionMethod(f(x)); endif;

end;
\end{verbatim}

\subsection{Polynomials of prime degree $p$ with Galois group $A_p$}

Let $f(x)$ be a polynomial in $\Q(t)$ as below $$f(x)= (n-1) x^n - n x^{n-1} + t.$$ The discriminant of $f(x)$
with respect to $x$ is
$$\D_f= (-1)^{\frac {n (n-1)} 2} n^n (n-1)^{n-1} t^{n-2} (t-1). $$
$\D_f$ is a complete square in $\Q$ if $(-1)^{\frac { (n-1)} 2} n t (t-1)$ is a complete square; see \cite{Se}
(pg. 44) for more on this family of polynomials. Let $n=23$. Then
$$\D_f= - 2^{22} \cdot 11^{22} \cdot 23^{23} \cdot t^{21} (t-1).$$
Hence, $\D_f$ is a complete square in $\Q$ if $G(t)= - 23 t (t-1)$ is a complete square. In other words, for all
those rational points on the curve $$y^2= G(t).$$ This is a genus 0 curve and can be parametrized as follows:
$$(y, t)= \left( - \frac {23m} {m^2+23} ,  \frac {23 } {(m^2+23} \right)$$
Consider $f(x)$ for $t =\frac {23 } {(m^2+23}$. Since we prefer to work with polynomials with integer coefficients
then take
$$f(x)= ( 22 m^2 + 506 ) x^{23} -( 23 m^2 + 529) x^{22} +23.$$
It is easily checked that $f(x)$ is irreducible over $\Q$ and its discriminant is
$$\D_f= 2^{22} \cdot 11^{22} \cdot 23^{44} \cdot m^2\, (23+m^2)^{22}$$ which is a complete square
in $\Q$. Thus, $Gal (f) $ is inside $A_{23}$. It is an simple calculus exercise to show that the number of real
roots of these polynomials is $\leq 3$. Hence, the Galois group is $A_{23}$.

We conclude with the following open problem:

\bigskip

\noindent \textbf{Problem:} \emph{Find a degree 23 polynomial $f(x) \in \Q[x]$ with exactly 7 real roots such that
$\D_f$ is a complete square in $\Q$ but $Gal (f) $ is not isomorphic to $A_{23}$.}

\medskip

\noindent The reader is probably aware that the above is an open problem of the inverse Galois problem, see
\cite{V}. Its solution would realize $M_{23}$ as a Galois group over $\Q$.

\section{Concluding remarks}

The algorithm suggested in this paper works very well with prime degree polynomials which have nonreal roots.
Since most polynomials have such roots this is an effective test to be implemented on all algorithms computing
Galois groups. Most computer algebra packages have already algorithms implemented to find the number of real roots
of a polynomial. For example in Maple the user can easily check using the command \verb"realroot( f(x) );". One
can also generalize this algorithm to any degree $n$ polynomial. However, in this case a more detailed analysis is
required.

\end{document}